\documentclass[11pt]{amsart}
\usepackage{geometry}                
\usepackage[parfill]{parskip}    
\usepackage{graphicx}
\usepackage{amssymb}
\usepackage{epstopdf}
\usepackage[all]{xy}
\title{Categories without structures}
\author{Andrei Rodin}

\begin {document}
\maketitle
\footnote{\emph{Acknowledgment}: I am very grateful to Sergei Artemov, Michael Detlefsen, Haim Gaifman, Daniel Isaacson, William Lawvere, Colin McLarty, Richard Pettigrew and Noson Yanofsky for their comments on earlier drafts of this paper and valuable discussions.}
\begin {abstract}
The popular view according to which Category theory provides a support for Mathematical Structuralism is erroneous. Category-theoretic foundations of mathematics require a different philosophy of mathematics. While structural mathematics studies Òinvariant formsÓ (Awodey) categorical mathematics studies covariant and contravariant transformations which, generally, donÕt have any invariants. In this paper I develop a non-structuralist interpretation of categorical mathematics and show its consequences for history of mathematics and mathematics education.
\end {abstract}
\section{Introduction: Renewing Foundations}
\paragraph{Foundations of mathematics, historically speaking, is the least stable part of this science. The Pythagorean theorem, for example, has stood firm since its early discovery while, during the same historical period, the foundations of geometry (and in particular the foundations of the Pythagorean theorem) have several times overturned and totally renewed.  Although these changes in foundations affected the theorem, there is still a sense in which the Pythagorean theorem always remained the same. (In the Conclusion of this paper I'll make this sense more precise). But the foundations of mathematics underlying different versions of this theorem did not remain the same. The traditional view that mathematics is about number and magnitude and the modern view that mathematics is about sets and structures differ radically. Actually saying that different foundations differ radically is pleonastic since one describes a difference as radical exactly when it concerns foundations rather than anything else. The above observation shows that the popular architectural metaphor of science, which describes science as an edifice with a solid foundation, is completely misleading when one talks about science from a long-term historical perspective. The renewal of foundations is not only compatible with the progress of science but also helps to make this progress possible.}

\paragraph{The notion of scientific progress assumes that once certain knowledge is acquired it later remains preserved and publicly available. It happens, of course, that certain beliefs, which at some point in history are generally seen as elements of current scientific knowledge, are later refuted and disqualified. However, the notion of scientific progress concerns knowledge itself, not our current beliefs about what does and what does not qualify as knowledge. Blurring the distinction between true knowledge and related beliefs would render the notion of progress incoherent. So let me now ignore the issue of belief revision and ask a different question:   Where is ready-made knowledge preserved and how does it endure through human history?}

\paragraph{According to Popper, \cite{Popper:1978} scientific knowledge and other products of the human intellect live in a special metaphysical domain that he calls the \emph{Third World}. The \emph{First World} on Popper's account is the world of physical processes and physical objects while the \emph{Second World} is that of mental states. Popper's rationale behind his notion of Third Worlds is that it avoids reducing knowledge to either mental states or physical processes:}

\begin{quote}
Knowledge in the objective sense consists not of thought processes but of thought \emph{contents}. It consists of the content of our \underline{linguistically}  formulated theories; of that content which can be, at least approximately, \underline{translated from one language into another} . The objective thought content is that which remains invariant in a reasonably good translation. Or more realistically put: the objective thought content is what the translator tries to keep invariant, even though he may at times find this task impossibly difficult. (Italic is Popper's, underlining mine)
\end{quote}

\paragraph{What is relevant to our present discussion here is not Popper's metaphysical argument but the way in which Popper thinks about thought contents in general and the content of scientific theories in particular. As a matter of course Popper doesn't identify content with its linguistic expression. He describes content as an invariant of linguistic translations of a given expression from one language into another. Using today's popular mathematical jargon we would say that Popper thinks here of the thought content as a linguistic pattern taken ``up to translation".}

\paragraph{I claim that Popper's notion of thought content fails to account for the long-term endurance of scientific, and in particular mathematical, knowledge. His theory better applies to the content of a religious doctrine rather than scientific content. A teacher of religion may indeed translate a sacred text of his religion to his less educated pupils doing his best to keep the original sense invariant. Even if the spirit of a religion is not generally reducible to the letter most developed religions use sacred texts as a means of preserving their identities over generations. But science proceeds very differently. A mathematical teacher - or at least a good mathematical teacher - doesn't try to transmit to her students the invariant content of some canonical text. Teaching the Pythagorean theorem today she doesn't ``try to keep invariant" what Euclid wrote about it some 2300 years ago but rather she relies on modern textbooks. If the notion of a canonical text makes a sense in science at all, it should be stressed that canonical scientific texts get quickly outdated, are revised, updated and periodically wholly rewritten.}

\footnote{Euclid's \emph{Elements} are often referred to as a typical example of a canonical mathematical text. It is often said that until recently people used this book as a Bible of geometry. In fact this alleged stickiness to Euclid's letter never existed. To see this it is sufficient to look more precisely into books titled Euclid's \emph{Elements} published before the 19th century. One finds a surprisingly diverse literature under this title. Early publishers and translators of Euclid's \emph{Elements} tried to produce a sound mathematical textbook rather than reproduce a canonical text. They didn't hesitate to improve on earlier editions of the Elements when they judged this appropriate. According to today's common standard the existing early editions of the \emph{Elements} don't qualify as different versions of the same text. Any of these people could get today a copyright as the author of his Euclid's \emph{Elements}. The notion of being an author is certainly changed since then. Today's canonical Greek edition of Euclid's \emph{Elements}  \cite{Euclides:1883-1886} was not produced  by Heiberg and his assistant Menge until the end of 19th century ; noticeably these people were philologists, not mathematicians. So the idea of reproducing Euclid's text literally and translate it into modern languages ``keeping its content invariant" is relatively recent; it is relevant to the history of mathematics rather than to mathematics itself. It should be also stressed that since the 17th century (a modernized version of) Euclid's \emph{Elements} was no longer the only geometry textbook on the market. Arnauld and some other authors produced their own original textbooks. }

\paragraph{Notice the difference between this latter kind of revision and the belief revision I mentioned earlier. Euclid's book in its original form is no longer in use in schools, but not because some of Euclid's propositions have been judged false by authorities of education. The fact that Euclid fails to meet today's standard of mathematical rigor is not the reason either, because elementary textbooks anyway do not meet and are not supposed to meet such a standard. The principal reason why Euclid's Elements is no longer used in schools is this: this book no longer provides a satisfactory basis for a study of more advanced and more specific branches of mathematics. It did this job perfectly  for quite a while but lost this capacity when mathematics essentially changed its shape. This dynamic is closely related to mathematical progress but it cannot itself be described as progress. Kids learning mathematics are hardly more intelligent today than they were a hundred years ago. Their learning capacities have hardly essentially increased. But today's kids should be prepared to use and further develop mathematics, which significantly progressed during the passed century. This is why they need a new curriculum.} 

\paragraph{This new curriculum cannot be just an extension of an older curriculum, because this would require an increase of pupils' learning capacities. So they should study a different mathematics from the beginning of their studies. This is why the evolution of mathematical curricula is not progressive in the precise sense of the term. }

\paragraph{True, older and newer elementary mathematics textbooks typically share some material. For example older and newer geometry textbooks usually include the Pythagorean theorem. Here is the Pythagorean theorem as it appears in Euclid's \emph{Elements}  (Proposition 1.47):} 
\begin{quote} \textbf{(1)} In right-angled triangles the square on the side subtending the right angle is equal to the squares on the sides containing the right angle.
\end{quote}
\paragraph{And here is how the same theorem appears in a modern textbook \cite{Doneddu:1965},( p.209, slightly modified):}
 \footnote{The original version reads \begin{quote} Two non-zero vectors x and y are orthogonal if and only if $(y - x)^2 = y^2 + x^2$ \end{quote} 
I use for my example only the ``only if" part of the statement. Euclid counts the converse of his 1.47 as a separate Proposition 1.48. Doneddu makes both propositions into one theorem. I have modified Donnedou's theorem because this obvious difference is irrelevant to my argument.}  

\begin{quote} \textbf{(2)}   If two non-zero vectors x and y are orthogonal then$(y - x)^2 = y^2 + x^2$ \end{quote}

\paragraph{In order to interpret the two propositions correctly, one must be careful. Euclid speaks here not about the areas of the squares but about the squares themselves: saying that the two smaller squares (taken together) are equal to the bigger square he means, roughly, that the bigger square can be composed out of pieces of the smaller squares.  The minus sign on the left side of  $(y - x)^2 = y^2 + x^2$ and the plus sign on the right don't stand for mutually inverse operations since the former operation applies to vectors while the latter applies to real numbers. Vectors, numbers, and operations with these things are construed here as structured sets. In order to interpret correctly not only the two statements but also their proofs much more needs to be said about their corresponding theories. In particular, much needs to be said about the structuralist set-theoretic foundations of mathematics developed in Doneddu's book and the foundations of Greek (and more specifically Euclid's) mathematics. But even without going into these details it is clear that the difference between \textbf{(1)}  and \textbf{(2)}  is not superficially linguistic. The two versions of the Pythagorean theorem differ in their foundations, i.e., differ radically.  And yet they express the same theorem! }

\paragraph {In the Conclusion of this paper I'll come back to this example and say more about it. But what has been already said already suffices for showing that Popper's account of thought content doesn't apply to mathematics as far as this science is viewed historically. Whatever the mathematical content might be, it cannot be described as an invariant of linguistic translation; the notion of linguistic translation doesn't help one to account for the long-term endurance of mathematical knowledge. The same is true of science in general. Unlike religious doctrines, poems, musical symphonies, and a few other inhabitants of Popper's Third World, scientific knowledge endures long-term through permanent revision rather than through mere repetition of linguistic patterns or retention of translational invariants, as Popper suggests. Above I described this revision as a pedagogical necessity. But it has also a philosophical aspect. I claim that the continual questioning, revision and renewal of foundations is a distinctive way in which science  endures through time and progresses. This ``unended quest" (to use Popper's word) concerns not only new yet unexplored domains of reality; it also concerns what is already known and well established.}

\paragraph{The renewal of foundations amounts to the dialectical refutation of older foundations and the dialectical positing of new foundations. This activity is philosophical rather than purely scientific. In this latter respect my view is traditional and qualifies as a form of foundationalism. But I also think that the notion of foundation does not make sense in abstraction from what it is (or what it is supposed to be) the foundation of. The historical performance of dialectically posited foundations crucially depends on what scientists (including mathematicians) do with foundations. So my scientific foundationalism implies the need for close cooperation between philosophy and science but definitely not the subordination of one of the two parties to the other. I subscribe to the following strong claim about the nature of scientific foundations:}

\begin{quote}A foundation makes explicit the essential general features, ingredients, and operations of a science, as well as its origins and generals laws of development. The purpose of making these explicit is to provide a guide to the learning, use, and further development of the science. A ``pure" foundation that forgets this purpose and pursues a speculative ``foundations" for its own sake is clearly a nonfoundation.\cite{Lawvere:2003}\end{quote}

\section{Claims and Slogans}
\paragraph{The above generalities help me to formulate the principal purpose of this paper. I shall consider Mathematical Structuralism as a particular project in foundations of mathematics but not as a speculative view on mathematics developed for its own sake. To proponents of Mathematical Structuralism in this latter sense I have nothing more specific to say except that I disagree with them on more general grounds, as explained in the above Introduction. I claim that Structuralism has been extremely successful in twentieth-century mathematics but now has been already worked out and should be abandoned. Recent developments in mathematics call (as ever) for renewal of foundations. My more specific claim concerns Category theory. This theory emerged in the mid-twentieth century within structural mathematics. Many proponents of Category theory believe that this theory is capable of providing better foundations for structural mathematics than the standard set-theoretic foundations used in the structural mathematics earlier (see Section 7 below).  I share with these people their enthusiasm about making Category theory into a new foundations of mathematics but I don't believe that mathematics built on these new foundations can or should be structural.}

 \footnote{By structural mathematics I understand mathematics built on foundations complying with principles of Mathematical Structuralism. A canonical example of structural mathematics is given by  \cite{Bourbaki:1939-1983}. }

 \paragraph{I claim, on the contrary, that Category theory brings about a new understanding of mathematics, which is connected to Structuralism historically and dialectically, but that is itself quite different. In what follows I shall try to present and develop this new view. I shall occasionally call this new view ``categorical" but leave it to others to invent a new \emph{ism} for it. However, following a fashion, I shall propose a polemical slogan. In the end of his \cite{Awodey:1996} Awodey puts forward the following structuralist slogan:}

 \begin{quote} The subject matter of pure mathematics is invariant form, not a universe of mathematical objects consisting of logical atoms. \end{quote} 
 
\paragraph{I counter this slogan with the following one:}

\paragraph{\emph{The subject matter of pure mathematics is transformation, not invariant form.}}

\paragraph{The rest of this paper is organized as follows. First, I briefly discuss Mathematical Structuralism, its historical origins and its relation to Set theory and Category theory. Here I explain reasons why MacLane, Awodey, and some other people believe that Category theory provides a support for Mathematical Structuralism. Then I provide my critical arguments against this view arguing that the notion of category should be viewed as generalization of that of structure rather than as a specific kind of structure. Further I analyze Lawvere's paper \cite{Lawvere:1966} on categorical foundations and show that the author begins this paper with a version of structuralist foundations but then proceeds with a very different project that still waits to be accomplished.  I conclude with an attempt to outline the new categorical view on mathematics explicitly.}

\section{Mathematical Structuralism}
\paragraph{Here is how a mathematical structure is described by a working mathematician for a philosophical reader:}

\begin{quote}All infinite cyclic groups are isomorphic, but this infinite group appears over and over again - in number theory, in ornaments, in crystallography, and in physics. Thus, the ``existence" of this group is really a many-splendored matter. An ontological analysis of things simply called ``mathematical objects" is likely to miss the real point of mathematical existence. \cite{MacLane:1996}\end{quote}

\paragraph{For a philosophical reader who doesn't know what the infinite cyclic group is (and also for the sake of my argument) I propose a modification of MacLane's example that amounts to a replacement of the words ``infinite cyclic group" by the words ``number three" and the word ``isomorphic" by the word ``equal":}

\paragraph{\emph{All threes are equal but this number appears over and over again - in number theory, in ornaments .... Thus the ``existence" of this number is really a many-splendored matter.}}

\paragraph{Indeed the familiar number three is just as promiscuous as the infinite cyclic group or perhaps even more promiscuous. It equally ``appears" (to use MacLane's word) both inside and outside mathematics: in a trio of apples, a trio of points, a trio of groups, a trio of numbers or a trio of anything else. As in MacLane's example, there is a systematic ambiguity between the plural and the singular forms of nouns in our talk about numbers. (Notice MacLane's talk about ``all infinite cyclic groups" and ``this infinite group" in the same sentence; in my paraphrase I talk similarly about a number.)  This analogy reveals a traditional aspect of Structuralism, which often remains unnoticed when people stress the novelty of this approach. Of course, this analogy doesn't allow for reduction of Structuralism to earlier views. But it allows one to see clearer what was really new in Structuralism. It was not the notion of ``many-splendored existence" stressed by MacLane in the above quote but a more specific notion of isomorphism, which plays in structural mathematics roughly the same role as the notion of equality (as distinguished from identity) plays in traditional mathematics.}

\paragraph{This point has been made clear by Hilbert in his often-quoted letter to Frege of December 29, 1899. Stressing the ``many-splendored" nature of structural theories (as we would call them today) Hilbert says:}

\begin{quote}[E]ach and every theory can always be applied to the infinite number of systems of basic elements. One merely has to apply a univocal and \underline{reversible} one-to-one transformation [to the elements of the given system] and stipulate that the axiom for the transformed things be correspondingly similar (quoted by \cite{Frege:1971}, underlining mine)\end{quote}

\paragraph{We see that Hilbert explicitly mentions here the reversibility condition, which implies that the given transformation is an isomorphism. Hellman  \cite{Hellman:forthcoming} quite rightly, in my view, recognizes Hilbert as a founder of Structuralism; however, in his official definition of Mathematical Structuralism, Hellman doesn't mention isomorphism explicitly: }

\begin{quote}Structuralism is a view about the subject matter of mathematics according to which what matters are structural relationships in abstraction from the intrinsic nature of the related objects. Mathematics is seen as the free exploration of structural possibilities, primarily through creative concept formation, postulation, and deduction. The items making up any particular system exemplifying the structure in question are of no importance; all that matters is that they satisfy certain general conditionsÑtypically spelled out in axioms defining the structure or structures of interest - characteristic of the branch of mathematics in question.\end{quote}

\paragraph{The reason why Hellman doesn't mention the notion of isomorphism here becomes clear from his remark concerning ``axioms defining the structure or structures [notice the plural - A.R.] of interest". Take axioms defining the notion of group for example. A group is any ``system" (to use Hellman's word) that consists of certain ``items" a , b, ... and binary operation $\oplus$ associating with every ordered pair of such items (possibly identical) a third item (possibly identical to one of those) from the same system such that the following axioms hold:}

\paragraph{\textbf{G1}: operation $\oplus$ is associative}
\paragraph{\textbf{G1}: there exists an item \emph{1} (called unit) such that for all \emph{a}  \emph{a}$\oplus$\emph{1} = \emph{1}$\oplus$\emph{a} = \emph{a} }
\paragraph{\textbf{G3}: for all \emph{a} there exists $a^{-1}$ (called inverse of  \emph{a}) such that $a\oplus a^{-1} = a^{-1}\oplus a$ =  \emph{1}}

\paragraph{These axioms are satisfied by many different groups. The infinite cyclic group mentioned above is just one example but there are many others.}

\footnote{We are now ready to spell out the precise definition: an infinite cyclic group is a group with an infinite number of elements and such that any of its elements is generated  by some distinguished element \emph{g} and its inverse $g^{-1}$. A group is said to be generated by a set of its distinguished elements called generators when every element of this group is a product of the generators. A canonical example of an infinite cyclic group is the additive group of whole numbers, which is generated by numbers 1 and -1.}   
  
\paragraph{Consider, for another example, a group of permutations of three letters \emph{A},\emph{B}, \emph{C} with the composition of permutations as group operation. This group comprises six different permutations (including the identical permutation). Since this latter group is finite it  cannot be isomorphic to the infinite cyclic group. So these two groups are different in the structural sense. This example shows that  the above axioms  \textbf{G1-3} describe structures of a particular type, not a particular structure.}

\paragraph{In order to give a sense to the expression ``type of structures" one needs to have the notion of structure beforehand. Axioms \textbf{G1-3}, or any other system of axioms determining some type of structure, cannot help one grasp the notion of structure unless one is already aware of the fundamental role of isomorphism. For the idea of a general description satisfied by different mathematical objects is obviously not unique to Structuralism; Euclid's axioms do the same job with respect to numbers and magnitudes. Stressing the higher importance of structures with respect to ``systems", the irrelevance of ``intrinsic nature" and relevance of ``structural relationships" cannot clarify the notion of mathematical structure by itself.} 

\section{Isomorphisms and ``Invariant Forms"}
      
\paragraph{A non-structuralist may observe that axioms \textbf{G1-3} are satisfied by a number of ``particular systems" (not structures so far!) called groups. Let now  \textbf{\emph{G}} be a class of such systems (i.e. groups), and consider the situation when some of these, say $G_1$ and $G_2$ are isomorphic. This actually means two things:} 

\paragraph{\textbf{I1}: elements of $G_1$  are in one-to-one correspondence with elements of $G_2$;}  

\paragraph{\textbf{I2}: for all elements $a_1$, $b_1$, $c_1$ from $G_1$ such that $a_1\oplus b_1 = c_1$ the corresponding elements $a_2$, $b_2$, $c_2$ from $G_2$  satisfy  $a_2 \otimes b_2 = c_2$  where  $\oplus$ is the group operation in $G_1$ and $\otimes$ is the group operation in $G_2$.} 

\paragraph{A one-to-one correspondence between elements of two given groups that satisfies \textbf{I2} is called (group) isomorphism. Groups are isomorphic if and only if there exists isomorphism between them. Notice that, given two isomorphic groups, there are, generally, many different isomorphisms between them. One should not confuse isomorphism as a particular correspondence between elements of two groups and isomorphism as an equivalence relation defined on some class of groups. Isomorphism in the latter sense holds between two given groups if and only if there exists an isomorphism in the former sense between these groups. As we can see, this terminology is slightly confusing but it is too common to try to change it.}  

\paragraph{Since  isomorphism is an equivalence relation it divides class  \textbf{\emph{G}} into sub-classes containing only isomorphic groups. One may ignore differences between isomorphic groups and get through this act of abstraction various notions of groups-qua-structures (not to be confused with the general notion of group as a type of structure!), in particular, the notion of infinite cyclic group. To facilitate the language and provide this reasoning with some intuitive support one may talk and think about any particular structure as a thing ``shared" by all members of the corresponding isomorphism class. On this basis one may claim that \emph{the items making up any particular system exemplifying the structure in question are of no importance} (as does Hellman in the above quote).  This claim describes the aforementioned abstraction, which can be called structural abstraction.  However, one cannot forget about these exemplifying systems completely because this would destroy the whole reasoning bringing about the notion of mathematical structure. Noticeably Hellman needs the auxiliary notion of system in order to describe what is a mathematical structure. One might think that this additional notion (no matter what one calls it) plays a role in philosophical talk about structural mathematics but plays no role in structural mathematics itself. In the next Section I shall argue that this is not the case.}

\paragraph{There is yet a different way of thinking about isomorphism (this will be already the third meaning of the term by our account!), which is common in current mathematical practice and particularly pertinent for categorical mathematics, as we shall later see.  One may think about a one-to-one correspondence between elements of groups $G_1$ and $G_2$, which satisfies condition \textbf{I2}, as a transformation  \emph{i}:  $G_1\rightarrow G_2$ of one group into another group. Since a one-to-one correspondence is a symmetric construction the choice of $G_1$ as the source and $G_2$ as the target of this transformation is in fact arbitrary. In other words one and the same isomorphism-\emph{qua}-correspondence gives rise to two isomorphisms- \emph{qua}-transformations  \emph{i}: $G_1\rightarrow G_2$  and \emph{j}: $G_1\rightarrow G_2$ , which run in opposite directions and cancel each other on both sides. This latter property means precisely the following:  the composition transformation $i\circ j$ resulting from the application of transformation\emph{j} after transformation \emph{i} sends every element of $G_1$ into itself and composition transformation $j\circ i$ sends every element of $G_2$ into itself (beware that none of the two conjuncts implies the other). Given these conditions each of transformations  \emph{i} and  \emph{j} is called the  \emph{inverse} of the other. Hence this definition: a transformation is called an isomorphism when it has an inverse. See Footnote 9 for a more precise categorical version of this definition.} 

\footnote{Notice that the order in which transformations are composed matters. I use here the so-called \emph{geometrical} notation where the composition is written in the ``direct" order. According to another notation called \emph{algebraic} the composition is written in the inverse order.}

\paragraph{Thinking about isomorphism as a reversible transformation allows one to think of a structure shared by given transformed systems as an ``invariant form", i.e. a form invariant under the given transformation. Hence Awodey's structuralist slogan quoted in the end of Section 2.  This slogan describes the structural abstraction in these alternative terms: only the invariant form matters, transformed systems don't.  As we shall see in what follows, the notion of isomorphism-\emph{qua}-transformation, which may seem redundant in the context introduced so far, becomes indispensable in categorical mathematics. Noticeably Hilbert in the above quote (Section 3) talks about isomorphism as transformation, not as a symmetric one-to-one correspondence.}

\section{Structures versus Abstract Objects; Collections versus Transformations}

\paragraph{Given an equivalence relation defined for a class of mathematical objects, Frege \cite{Frege:1884} considered the possibility of replacing each obtained equivalence class by a single object through an act of abstraction.}

\footnote{Frege's example is the concept of direction built, as follows. One considers the class of all straight lines on a Euclidean plane and the equivalence relation ``is parallel". Then one associates a single abstract  concept called \emph{direction} with each isomorphism class of parallel lines.} 

\paragraph{Frege calls the result of this procedure an abstract object, not a structure, and indeed he doesn't think about this outcome as a structure. So we need a further effort for distinguishing structural  abstraction from other types of mathematical abstraction. To this end, let us first consider this question: What are elements of a group- \emph{qua}-structure? For the reason that I have already explained we don't want these elements to have anything like an ``intrinsic nature". So they should be just ``items" or ``abstract elements"; the predicate ``abstract" refers here to the act of abstraction through which the notion of group- \emph{qua}-structure is obtained. However, we still need to make some assumptions about these things. We want them to be many and form (or belong to) well-distinguishable collections. Since we want to use the same notion of collection for different purposes we don't want the collected elements to be related in a specific way. This will give us the freedom to stipulate any relation between elements by fiat using the same notion of collection.}

\paragraph{This is an important point where Structuralism meets Set theory. Having a notion of set at our disposal we are in a position to give the standard structural definition of group as a ``structured set", namely a set provided with a binary operation satisfying axioms  \textbf{G1-3} given above. There is a standard way to account for algebraic operations as relations that I shall not explain here.}

\paragraph{As we have seen, the notion of isomorphism plays a crucial role in structural abstraction, which brings about new mathematical objects (namely, new mathematical structures). Importantly isomorphisms do not disappear when a given act of structural abstraction is accomplished and a new mathematical structure is well-defined. Mathematicians think about abstract groups and other abstract structures as given in an indefinite number of isomorphic copies, not as unique objects. As I have already stressed, people think similarly about numbers in traditional arithmetic (see Section 3). This, in my view, is the principal point where Frege's notion of abstraction fails to account for structural abstraction as this latter notion has been developed in twentieth-century mathematics. Reasoning ``up to isomorphism" doesn't amount to the strict identification of isomorphic structures; it rather amounts by replacement of traditional equality by isomorphism in appropriate contexts. From a mathematical (as distinguished from logical and philosophical) viewpoint the question whether or not two structures are identical is just as pointless as the question whether or not two equal numbers are identical. A sound mathematical question about two given numbers is whether or not they are equal. A sound mathematical question about two given structures is whether or not they are isomorphic.}

\footnote{I elaborate on this issue in \cite{Rodin:2007}.}

\paragraph{Set theory makes the talk of isomorphism as transformation redundant because the notion of one-one correspondence may be analyzed set-theoretically in terms of pairs of elements. However in many important  mathematical contexts the notion of transformation is widely used anyway: groups of (reversible) transformations are abundant and geometry and also in physics. As far as foundations of mathematics are concerned we have an important choice here: either to (i) consider the notion of collection as more fundamental than that of transformation and reduce the latter to the former or to (ii) consider the notion of transformation as more fundamental and reconstruct the notion of collection on this basis. The former option brings (some version of) set-theoretic foundations of mathematics. The idea of categorical (i.e. category-theoretic) foundations amounts to taking the latter option. However the project turns to be non-viable unless one takes into the account other transformations than isomorphisms.}

\section{Homomorphisms}

\paragraph{Given a type of structures it is always possible to define a general notion of map between structures of the given type. I shall discuss first the case of general maps between groups; such maps are called \emph{homomorphisms} or more precisely \emph{group homomorphisms}.  Then I shall say few words about general maps between structures of different types. The term ``homomorphism" is traditionally reserved for groups (apparently because this case was studied first), although, as its etymology suggests, it could also be used for structures of different types like the term ``isomorphism". So in what follows I shall use the term ``homomorphism" in the sense of general map between structures of some given type.}

\paragraph{The notion of group homomorphism generalizes upon that of group isomorphism in the following way: instead of one-to-one correspondence between elements of groups $G_1$, $G_2$, one considers a more general kind of correspondence that is allowed to be many-to-one (but not one-to-many). In other words, one considers a \emph{function} (in the set-theoretic sense of the term) \emph{f}: $S_1 \rightarrow S_2$ from the set $S_1$ of elements of $G_1$ to the set $S_2$ of elements of $G_2$. Condition \textbf{I2} from Section 4 remains the same; notice that it can be satisfied when elements $a_1$, $b_1$ are different but elements $a_2$, $b_2$ are the same.}

\paragraph{Group homomorphism and similar general maps between structures of other types are colloquially called ``structure preserving". This is somewhat misleading because if such maps preserve anything at all it is a type of structure but not a particular structure. Think about this trivial example: for all groups $G_1$, $G_2$ there exist a homomorphism \emph{h}: $G_1 \rightarrow G_2$  which sends every element of $G_1$ to the unit of $G_2$.  This homomorphism ``destroys all information" about $G_1$ reducing its image to a single element;  it doesn't provide any information about $G_2$ either.}

\paragraph{Actually the example of group homomorphism doesn't straightforwardly generalize  to maps betweens structures of different types. For given a type of structure there are, generally, different ways to define maps between structures of the given type (some of which may be reasonable and some other not). Such maps can be of different kinds. Usual maps between topological spaces, i.e., general continuous transformations, do not  \emph{preserve} topological structure (in the same sense in which group homomorphisms are said to preserve group structure) but  \emph{reflect} it:  the \emph{inverse} image of any open set under a given continuous transformationis is always open while the direct image of an open set can be closed. In the case of isomorphism the difference between reflection and preservation of structural properties disappears. This fact shows that thinking about homomorphisms as ``imperfect isomorphisms" can be misleading; at the very least one should not forget that a given structural isomorphism may ``loose its perfection" in two different ways.}    

\paragraph{I shall now argue that homomorphisms, generally, don't allow for invariants in anything like the same sense in which isomorphisms do so. Let us try to replace isomorphisms by homomorphisms in the process of structural abstraction described in Section 4 and see what happens. One might expect to get in this way a generalized notion of structure but this doesn't work. Recall the first step: given class \textbf{\emph{G}} of groups we have divided it into equivalence subclasses of isomorphic groups. Two groups are isomorphic if and only if there exists isomorphism (i.e., a reversible transformation) between them; clearly this is an equivalence relation. Let me (for the sake of argument) call two groups \emph{homomorphic} if and only if there is a homomorphism between them. Although this latter relation is also an equivalence, one can see the difference: since \emph{all} groups are homomorphic (see the above example of group homomorphism) one cannot use this equivalence for dividing \textbf{\emph{G}} into equivalence subclasses! Saying that two given groups are homomorphic is tantamount to saying that the given groups are groups. So the relation of homomorphism just introduced (not to be confused with the standard notion of homomorphism as transformation) doesn't make sense.}  

\paragraph{In order to see the reason for this failure, note that the existence of homomorphism of the form $G_1 \rightarrow G_2$ doesn't imply the existence of homomorphism of the form $G_2 \rightarrow G_1$. This means that in the case of homomorphism (unlike that of isomorphism) the difference between the source and the target of the given transformation matters. But the relation of homomorphism tentatively introduced above doesn't take this difference into account. It forgets the difference between isomorphic and non-isomorphic groups and thus confuses their structural properties and offers no replacement.} 

\paragraph{A more reasonable choice of relation associated with a given homomorphism  \emph{h}: $G_1 \rightarrow G_2$ would be that of non-symmetrical relation $>$ such that $G_1 > G_2$ holds just in case there is a homomorphism of the form \emph{h}: $G_1 \rightarrow G_2$. However, since $>$ is asymmetric it is not an equivalence and so doesn't allow one to proceed further with the structural abstraction or anything similar.}

\paragraph{We see that homomorphisms cannot do the same job as isomorphisms: the reversibility condition stressed by Hilbert in the above quote (Section 3) turns out to be crucial for structural abstraction. One cannot reason ``up to homomorphism" in anything like the same way in which people reason up to isomorphism doing structural mathematics. Since ``invariant" in the given context is just another word for structure it is clear that homomorphisms, generally, don't have invariants in anything like the same sense in which isomorphisms and groups of isomorphisms do so.}

\footnote{I mention here groups of isomorphisms (not to be confused with isomorphisms of groups!) because they are very important in geometry and physics. I mean  groups of geometrical transformations of a given space. Only \emph{reversible} geometrical transformations, i.e. geometrical isomorphisms, of a single object (the given space) form groups (with the composition of transformations as group operation) because in this case the reversibility is equivalent to the existence of inverse elements. So the talk of invariants of groups, which is so important for structural approaches in physics, concerns only reversible transformations and doesn't apply to geometrical (or other) transformations in general. A non-mathematical reader may skip the reference to groups of isomorphisms in this part of the paper. I shall explain the idea of group of isomorphisms more clearly in categorical terms in Section 8.}

\section{Structuralists Motivations behind Category Theory}
\paragraph{The emergence of Category theory in he 1940s and its further development in the context of structural mathematics was related to a growing awareness of the role of general maps (not only isomorphisms). I shall not explain here the precise mathematical context in which this theory first proved useful but only mention that the notion of category generalizes upon such examples as the class of all sets and all functions, all groups and all group homomorphisms, all topological spaces and all continuous maps (not only reversible ones!) between topological spaces. This is a simple theorem \cite{MacLane:1996} that a class of structures of any fixed type provided with an appropriate notion of general map form a category. Generally, a category comprises a class of objects and a class of composable maps (called in Category theory \emph{morphisms}) for every ordered pair of objects, which are subject to few natural axioms. Given two different categories one defines a notion of map between categories. Such maps are called \emph{functors};  the usual definition of functor is based on the same idea as the definition of group homomorphism given in the previous Section: a functor sends each object of the source category into an object of the target category and each morphism of the source category into a morphism of the target category in such a way that composition of morphisms is ``preserved" in the same sense in which the group operation is said to be preserved by a group homomorphism. Using the notion of functor one may consider various categories of categories, i.e. categories such that their objects are themselves categories. One may also consider categories objects of which are functors. The above standard description of basic categorical concepts is structuralist in spirit. In Section 11 I shall describe functors and categories anew from a foundational and ``more categorical" viewpoint.}

\paragraph{The idea of categorical foundations as viewed from a structuralist perspective amounts to recovering all the relevant properties of any structure of any given type through properties of the category of (all) structures of this given type. In the case of the category of sets this provides an alternative (category-theoretic) Set theory: one first conceives of sets as abstract objects and stipulates that they form a category; then one stipulates desired properties of this category, which make this category ``into" the intended category of sets. This result (see \cite{Lawvere:1964}) shows that a reasonable notion of collection (set) can be developed on the basis of that of transformation (morphism of sets) but not only the other way round.} 

\paragraph{The growing popularity of Category theory as a common (albeit certainly not unique) ``language" of contemporary mathematics as well as the continuing efforts of building categorical foundations of mathematics are generally seen as a further step of the structuralist project briefly described above. I agree with this view so far as it does not require preserving the basic principles of Mathematical Structuralism (as specified above) in the new categorical setting. In my understanding, these developments diverge from Mathematical Structuralism and tend towards a very different view on mathematics and science in general. Before I describe this new view, let me explain reasons why categorical foundations appear to many as a version of structural foundations. In the next Section I shall show that this impression is wrong.}

\paragraph{As I have explained in Section 5, the notion of set plays a special role in structural mathematics. This explains why Set theory itself is rarely seen as a structural theory on equal footing with, say, Group theory. As Hellman \cite{Hellman:forthcoming} rightly remarks:}

\begin{quote}[D]espite the multiplicity of set theories (differing over axioms such as wellfoundedness, choice, large cardinals, constructibility, and others), the axioms are standardly read as assertions of truths about ``the real world of sets" rather than receiving a structuralist treatment.\end{quote}

\paragraph{The structural notion of group explained above is usually construed as a ``set with a structure" or ``structured set" rather than a pure structure (whatever this might mean); the underlying set of a given group is thought of as a background supporting the structure rather than a part of this structure. This way of thinking in mathematics is reminiscent of Aristotle's metaphysics of Matter and Form. The need for the set-theoretic Matter to do structural mathematics becomes clear from our analysis given in Section 5, but the presence of this ingredient doesn't fully comply with the philosophy of Mathematical Structuralism, which purports to make mathematical objects into pure forms (structures) and leave anything like their ``background" outside mathematics. The desired ``purely structural" mathematics would deal only with the ``invariant Form" and require no set-theoretic Matter. Historical evidence of such an attitude can be found in what Dieudonne (under the name of Bourbaki) says in his structuralist manifesto \cite{Bourbaki:1950} about set-theoretic difficulties:}

\begin{quote}The difficulties did not disappear until the notion of set itself disappears ... in the light of the recent work on the logical formalism. From this new point of view mathematical structures become, properly speaking, the only ``objects" of mathematics.\end{quote}

\paragraph{I don't believe that Dieudonne's claim concerning the alleged ``disappearance" of sets is justified but the quote clearly demonstrates such an intention.}

\paragraph{In this context the idea of accounting for relevant properties of mathematical structures only in terms of structure-preserving maps between these structures independently of any set-theoretic background, i.e., the idea of categorical foundations, indeed may look like a further step in the structuralist direction. Hence the popular view according to which categorical mathematics is the desired purely structural mathematics.} 

\paragraph{Remarkably, Category theory did never make it into Bourbaki's \emph{Elements} \cite{Bourbaki:1939-1983}, which is the most systematic attempt to develop structural mathematics ever undertaken. This is in spite of the fact that both founders of Category theory, Eilenberg and MacLane, were eventually involved in the Bourbaki group, so all the principal members of this group were well aware about their work.  This fact is often seen as a historical puzzle but in my view it is not. For, as we shall shortly see, categorical foundations of mathematic are not and cannot be anything like the structural foundations developed by Bourbaki in his fundamental work. So in order to incorporate Category theory into their \emph{Elements} Bourbaki would need to abandon his basic structuralist principles and engage  himself into a very different foundational project.}  

\paragraph{One may agree that Bourbaki's version of Structuralism is incompatible with categorical foundations of mathematics but argue that some other variety of Structuralism is appropriate for building such foundations. The force of this argument evidently depends on what exactly one understand by Mathematical Structuralism.  There are several varieties of this notion discussed in the recent literature \cite{Hellman:2001}. Rather than elaborate on such varieties I tried to explicate in this paper a general notion of Structuralism that would reflect the structural character of mathematics 
of 20th century mathematics broadly conceived. My arguments concerning Structuralism and Category theory apply to this general notion of Structuralism but not only to simple Bourbaki-like examples of structures and maps between structures, which are given here for a non-mathematical reader.} 
\paragraph{One may argue that my notion of Structuralism is nevertheless too restrictive and doesn't reflect the structural character of modern mathematics in full. Even if in this case the issue may look merely terminological I would stress the need to define one's general notions of structure and Structuralism explicitely and precisely. Distinguishing between multiple varieties of Structuralism doesn't help one to meet this requirement unless one addresses the question What these different varieties are varieties of? What I want to stress in this paper is a conceptual difference between the ``classical'' structuralist thinking exemplified by \cite{Hilbert:1980} and \cite{Bourbaki:1939-1983}, on the one hand, and some developments in Category theory, on the other hand. Leaving terms ``structure'' and ``Structuralism'' without any precise definition and using them in the broadest possible sense  can  hardly be helpful for showing such a difference. If Category theory indeed brings about a new philosophy of mathematics this new philosophy needs a new vocabulary.}

\section{Categories versus Structures; Embodiement of Mathematical Concepts}

\paragraph{Categories of structures like the category of groups, topological spaces, etc. capture the notion of type of structure, not the notion of singular structure. Particular structures (identified up to isomorphism) may be often also rendered as categories but in this case their morphisms are no longer structure-preserving maps. For example, a particular group (like the infinite cyclic group mentioned above) can be presented as a category with just one object such that all of its morphisms (going from this object to itself) are isomorphisms. The group operation is given by composition of morphisms; the existence of unit follows from the definition of a (general) category and the existence of inverse elements follows from the fact that all morphisms of the given category are reversible.}

\footnote{Categorical definition of isomorphism resembles the definition of reversible transformation given in the end of Section 4. However, it doesn't involve a reference to elements. Think about groups  $G_1$,  $G_2$ as objects of a category and modify the definition of Section 4 in this way: $i \circ j = 1_{G_1}$ and $j \circ i = 1_{G_2}$ where  $1_{G_1}$ is the identity morphism of $G_1$ and $1_{G_2}$  is the identity morphism of $G_2$. The rest of the definition remains the same.} 

\paragraph{This simple example shows that categorical morphisms can but should not be structure-preserving maps. Moreover, the above categorical presentation of group, unlike its standard set-theoretic presentation, is not structuralist in character. For the standard structuralist presentation involves this idea: an abstract group can be ``exemplified" by what Hellman calls ``particular systems", like systems of numbers, systems of geometrical motions and so on and so forth. Of course, when one pictures elements of a given group as loops rather than dots this does not produce any conceptual change by itself. But given the above categorical presentation of a group, and using standard category-theoretic means, one can do something other than keep saying that morphisms of the given category (i.e., the given group) stand or may stand for something else than themselves. Namely, one may consider functors from the given group-category into some other categories, which in their turn present (rings or fields of) numbers, geometrical spaces, etc. This provides a much more precise idea of ``standing for" in each particular case than the general structuralist rhetoric. In the structuralist setting the notion of exemplification remains meta-theoretical and escapes a precise mathematical treatment. But in the categorical setting this notion becomes a proper part of the given mathematical construction. Instead of saying that $A$ stands for $B$ one considers functors of the form $A\rightarrow B$ and treats these functors on equal footing with ``internal" morphisms of $A$ and $B$.}

\footnote{A further step of such categorical analysis amounts to considering the full category of functors of the given form; such a functor category provides a precise information about how $A$ translates into $B$.}

\paragraph{In my understanding, this latter type of mathematical thinking has little if anything to do with structural abstraction. A principal epistemic strategy of Structuralism is to capture what various `particular systems" share in common, namely their ``shared structure". The corresponding categorical strategy can be described in this way: look how ``particular systems" translate into each other. Unlike the structuralist strategy this categorical strategy doesn't make the particular systems less important. Given morphism $A\rightarrow B$ there is, generally, no reason to think of $A$ and $B$  ``up to" some equivalence and dispense with $A$ and $B$  in favor of their shared structure or anything else.  As I have already shown in Section 6 the notion of thinking ``up to homomorphism" is plainly unsound.} 

\paragraph{Let us now consider the case when a category presents a type of structure rather than a singular structure. To analyze this case I shall use the notion of \emph{embodiment}, which I have introduced elsewhere \cite{Rodin:forthcoming}. As we have seen in Section 5 a mathematical structure cannot be identified with its corresponding abstract concept: something else is needed in order to make a given concept into a mathematical object. Kant would call this additional element an intuitive construction; I use the word ``embodiment" for a similar purpose but in a different mathematical context. We have seen how the notion of structure allows for making a concept describing different ``particular systems" into a single mathematical object (single up to isomorphism, of course). As we have seen in Section 6 this structuralist method of embodiment doesn't work for \emph{types} of structure. While the concept ``infinite cyclic group" can be embodied into a single structure, the concept ``group" cannot; ``the group" unlike ``the infinite cyclic group" is not a name of unique (up to isomorphism or otherwise) mathematical object. However the \emph{category} of (all) groups is a single mathematical object like number 3, the infinite cyclic group or, say, the Euclidean plane. Each of these objects has a \emph{many-splendored existence} (to use MacLane's word), so its singleness must be understood appropriately. But I want now to stress a different point: the way in which all isomorphic cyclic groups are made into a single object with the notion of structure and the way in which all groups are made into a single object with the notion of category are essentially different. While the former involves structural abstraction the latter involves a different kind of abstraction, which I shall call \emph{categorical}. Roughly, categorical abstraction amounts to the following: one forgets about the fact that groups have elements and consider only how they map to (i.e. transform into) each other with appropriate morphisms; a relevant notion of element is recovered in this categorical setting only later on. Obviously the two kinds of abstraction are quite different. I shall say more about categorical abstraction in the Conclusion.}

\paragraph{A category in which morphisms (including identity morphisms) form a set (in the technical sense of the term) is called \emph{small}. Small categories can be thought of as structures on their own. The corresponding type of structures is defined straightforwardly: one takes a set of elements called morphisms, stipulates appropriate primitive relations between elements of this set, and spells out the necessary axioms (see the next Section for more details). Thus small categories like groups can be thought of as structures of a specific type. Noticeably, this straightforward approach doesn't work in the case of large categories corresponding to types of structures - think again of the category of groups or the category of all small categories. Since morphisms of such categories form proper classes they cannot be described as structured sets. Although this may look like a minor technical difficulty, which can be resolved by an appropriate generalization of the usual notion of structure, this difficulty provides additional evidence that the structural approach, generally, doesn't work in Category theory. Instead of thinking of categories as structures (or generalized structures) of a particular type, it seems to me more reasonable to reverse the order of ideas and think of structures as categories or categorial constructions of a particular type. An immediate suggestion would be to identify structures with small categories. A more elaborate suggestion by Lawvere (in person) is to identify a structure with a functor from a small category to a large ``background" category, say, that of sets.}  

\paragraph{To conclude this Section, let me stress that categories don't always represent particular structures or particular types of structure. Examples of this latter kind are today so popular only because they connect the new categorical mathematics with the older structuralist mathematics. But categorical mathematics also involves concepts and constructions that were first developed in a categorical setting, for example that of Grothendieck topology. One may expect that the further development of categorical mathematics will make such ``purely categorical" concepts better known and more useful in various branches of mathematics; then the link between the categorical mathematics and its structural predecessor will become a historical nd philosophical rather than mathematical issue.}

\section{``The category of categories"}
\paragraph{The idea of categorical foundations amounts to taking the notions of category, functor and/or some other related categorical notions as primitive and recovering the rest of mathematics on this basis. What are possible ways of realizing this project? In which precise sense can one consider category-theoretic notions as primitive? A way to do this, which immediately suggests itself, is to use in categorical foundations a modern version of Hilbert-style axiomatic method after the example of standard set-theoretic foundations. }

\paragraph{Consider a class of things called morphisms and three primitive relations: one that associates with every given morphism its source, one that associates with every given morphism its target, and, finally, one that associates with some (ordered) pairs of morphisms a third morphism called the composition of the given two morphisms. Then we need axioms to ensure that sources and targets of morphisms behave as identity morphisms (i.e. as \emph{objects}), that two given morphisms are composable if and only if the target of the first morphism coincides with the source of the second morphism, and some other similar axioms. Finally we should assume that the composition of morphisms is associative. For the full list of such axioms I refer the reader to \cite{Lawvere:1966}. The axiomatic theory  just described this author calls the \emph{elementary theory of abstract categories}.}  

\paragraph{Lawvere's paper begins as follows:}
\begin{quote} In the mathematical development of recent decades one sees clearly the rise of the conviction that the relevant properties of mathematical objects are those which can be stated in terms of their abstract structure rather than in terms of the elements which the objects were thought to be made of. The question thus naturally arises whether one can give a foundation for mathematics which expresses wholeheartedly this conviction concerning what mathematics is about, and in particular in which classes and membership in classes do not play any role. \end{quote}

\paragraph{We see that Lawvere embraces Mathematical Structuralism here but at the same time rejects set-theoretic (and even more general class-based) foundations of mathematics. Since the Hilbert-style axiomatic method is essentially structural (see Section 3 above) Lawvere's method of building his \emph{elementary theory of abstract categories} perfectly  fits his stated purpose. After the introduction of the axioms of the \emph{elementary theory} and providing some definitions on their basis Lawvere says:}

\begin{quote}By a category we of course understand (intuitively) any structure which is an interpretation of the elementary theory of abstract categories, and by a functor we understand (intuitively) any triple consisting of two categories and a rule $T$ which assigns, to each morphism $x$ of the first category, a unique morphism $xT$ of the second category in such a way that ... \end{quote}

\paragraph{(follow the conditions of being structure-preserving). A problematic aspect of this first part of the paper concerns Mayberry's argument that Lawvere's   \emph{elementary theory} like any other theory built with the Hilbert-style axiomatic method requires some primitive (non-axiomatic) notion of collection, which cannot be identified with that of category  \cite{Mayberry:2000}. The argument implies that the  \emph{elementary theory} and the corresponding elementary notion of category cannot be a genuine foundation. I agree with Mayberry on this point (this follows from my understanding of the relationships between Structuralism and Set theory explained in the beginning of Section 5), but unlike Mayberry I think that such a primitive notion of collection is dispensable in foundations of mathematics along with the Hilbert-style structural axiomatic method itself. In what follows I shall sketch a different version of axiomatic method that seems to me more appropriate for categorical foundations. Let me now return to Lawvere  \cite{Lawvere:1966}.}

\paragraph{Lawvere's \emph{elementary theory} is a preparatory step towards another theory of categories, which Lawvere calls \emph{basic theory}. My claim is that unlike the elementary theory the \emph{basic theory} is not structural, at least not in a similar sense. If I am right this shows that the main content of \cite{Lawvere:1966} in fact doesn't agree with the structuralist agenda announced by the author in the beginning of his paper: Lawvere begins with structural reasoning but then proceeds with a very different agenda, which can be described as genuinely categorical.}

\paragraph{The \emph{basic theory} begins with a re-introduction of the notion of functor:}
\begin{quote}Of course, now that we are in the category of categories, the things denoted by the capitals will be called categories rather than objects, and we shall speak of functors rather than morphisms.\end{quote}

\paragraph{This may sound like a mere terminological convention (rather than an alternative definition) but in fact it signifies a sharp change of perspective. The idea is now the following: given a preliminary notion of category (through the \emph{elementary theory}), conceive of category \textbf{\emph{C}}  of \emph{all} categories; then pick up from \textbf{\emph{C}}  an arbitrary object $A$ (i.e., an arbitrary category) and finally specify $A$ as a category by internal means of \textbf{\emph{C}},  stipulating additional properties of \textbf{\emph{C}}  when needed. More precisely it goes as follows (I omit details and streamline the argument). Stipulate the existence of terminal object  \emph{1} in \textbf{\emph{C}}, i.e., the object with exactly one incoming functor from each object of \textbf{\emph{C}}. Then identify objects (= identity functors) of $A$ as functors in \textbf{\emph{C}} of the form \emph{1}$\rightarrow A$. Stipulate also the existence of initial object \emph{0}, i.e. the object with exactly one outgoing functor into each object of  \textbf{\emph{C}}.  Then consider in  \textbf{\emph{C}} object  \emph{2} of the form \emph{0}$\rightarrow$\emph{1}  and stipulate for it some additional properties among which is the following: \emph{2} is a universal generator which means that:}

\paragraph{\textbf{G} (generator):  for all $f$, $g$ of the form:} 
$$\def\dar[#1]{\ar@<2pt>[#1]^f \ar@<-2pt>[#1]_g} \xymatrix{ A \dar[r] & B }$$
\paragraph{and such that $f \neq g$ there exist $x$ such that:}     
$$\def\dar[#1]{\ar@<2pt>[#1]^f\ar@<-2pt>[#1]_g} \xymatrix{ \emph{2}\ar[r]^x & A \dar[r] & B }$$

\paragraph{and $xf \neq xg$.}

\paragraph{\textbf{U} (universal): if any other category $N$ has the same property, then there are $y$, $z$ such that:}

$$\xymatrix{ A \ar@<2pt>[r]^y & B \ar@<2pt>[l]^z }$$

\paragraph{and  $yz = \emph{2}$.} 

\paragraph{This allows Lawvere to identify functors (morphisms) of $A$ as functors of the form \emph{2}$\rightarrow A$ in $C$. The fact that \emph{2} is the universal generator (it is unique up to isomorphism as follows from the above definition) assures that categories are determined ``arrow-wise": two categories coincide if and only if they coincide on all their arrows. This new definition of functor also allows one to make sense of the notion of a component of a given functor of the form $h$: $A\rightarrow B$ , which in the elementary theory is understood as a map $m$ sending a particular morphism $f$ of $A$ into a particular morphism $g$ of $B$ . In the basic theory, $m$ turns into this commutative  triangle:}

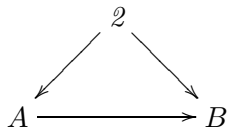
\footnote {A categorical diagram is said to commute or be commutative when the compositions of all morphisms shown at the given diagram produce other morphisms shown at the same diagram in appropriate places, so that any ambiguity about results of the compositions is avoided. For example, saying this triangle 

$$\xymatrix{&B\ar[dr]^g\\ A\ar[ur]^f\ar[rr]_h && C}$$

is commutative is simply tantamount to saying that $fg = h$. Morphisms resulting from composition of shown morphisms can be omitted at a commutative diagram when this doesn't lead to an ambiguity. For example, saying this square

$$\xymatrix{A \ar[r]^g & B \\ C \ar[r]_h \ar[u]^f & D\ar[u]_i}$$

is commutative is tantamount to saying that $fg = hi$.}

$$\xymatrix{&\emph{2}\ar[dl]\ar[dr]\\ A\ar[rr] && B}$$ 

\paragraph{This, once again, significantly changes the whole perspective: categories and functors are no longer built ``from their elements" but rather ``split into" their elements when appropriate. Although the notion of functor as a structure-preserving map can be recovered in this new context it no longer serves for defining the very notion of functor. Rule $T$ used by Lawvere for defining functors in the \emph{elementary theory} disappears in the \emph{basic theory} without leaving any trace.}

\paragraph{Further consider this triangle which Lawvere denotes \emph{3}:} 
  
$$
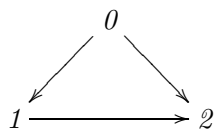
\xymatrix{&\emph{0}\ar[dr]\ar[dl]\\ \emph{1}\ar[rr] &&\emph{2}}$$

\paragraph{(It should satisfy a universal property which I omit). \emph{3} serves for defining composition of morphisms in our ``test-category" $A$ as a functor of the form \emph{3}$\rightarrow A$ in $C$. Finally, in order to assure the associativity of the composition Lawvere introduces category \emph{4}, which is pictured as follows:}

$$\xymatrix{&\emph{3}&\\ \emph{0}\ar[ur]\ar[rr]\ar[dr] &&\emph{2}\ar[ul]\\&\emph{1}\ar[ur]\ar[uu]&}$$ 

\paragraph{(The associativity concerns here the path  $\emph{0}\rightarrow\emph{1}\rightarrow\emph{2}\rightarrow\emph{3}$.)} 
\paragraph{This construction provided with appropriate axioms makes A into an ``internal model" of the elementary theory in the following precise sense: If $F$ is any theorem of the \emph{elementary theory}, then ``for all $A$, $A$ satisfies $F$" is a theorem of the \emph{basic theory}.}

\footnote{Isbell in his review \cite{Isbell:1967} of Lawvere's \cite{Lawvere:1966} points to a technical flaw in Lawvere's proof of this theorem . This flaw is fixed, in particular, in \cite{McLarty:1991}.}

\paragraph{The following analogy with the set-theoretic mathematics helps to clarify the role of categories of categories in foundations. As long as the notion of set is not supposed to provide a foundation for mathematics, one thinks of sets after examples of sets of numbers, sets of points, and the like. But in a foundational axiomatic theory of sets like ZF there are no other sets but sets of sets, and every mathematical object like a number or a point is supposed to be a set. Similarly in a foundational axiomatic theory of categories there are no other categories but categories of categories and every mathematical object is ultimately a category.}

\section{Functorial Semantics, Sketch Theory and Internal Language}
\paragraph{In order to see that Lawvere's basic theory unlike his elementary theory is not based on structuralist principles, and then to get an idea of non-structuralist principles behind this theory, it is instructive to take into consideration two similar approaches:  \emph{Functorial semantics} developed by the same author elsewhere \cite{Lawvere:1963-2004} and  \emph{Sketch theory} founded by Ch. Ehresmann in the 1960s and later developed by other people (see  \cite{Wells:1993} for an overview and further references).}

\paragraph{Functorial semantics involves the presentation of mathematical theories as categories of a special sort;  models of a given theory are functors from the theory to the background category of sets or another appropriate topos. The very idea of ``interpretation" or ``realization" of a given theory in a set-theoretic background obviously comes from the standard (due to Tarski) Model theory.  Lawvere's functorial semantic can be seen as a category-theoretic version of the same basic construction. However, as we shall now see, this technical update comes with a significant revision of the structuralist background of Tarski's Model theory inherited from Hilbert's notion of axiomatic method.}

\paragraph{In order to determine a theoretical structure, an axiomatic theory should  be \emph{categorical} , i.e., to have models that are all isomorphic. (Beware that this older sense of the term ``categorical" has nothing to do with Category theory!). True, not all axiomatic theories built by the standard method satisfy this requirement; also true, non-categorical theories are usually not disqualified solely on this basis. Anyway, in the standard setting the categoricity of axiomatic theory is commonly (and usually as a matter of course) viewed as an epistemic gain while the lack of categoricity is viewed as a problem. As long as one commits oneself to Structuralism such an attitude is understandable:  when a set of axioms fails to specify a model up to isomorphism it fails to specify a structure. Saying that a non-categorical theory determines many structures rather than one structure is somewhat misleading because such a theory, strictly speaking, doesn't specify any structure at all (cf. Section 3 above).}  

\paragraph{In the case of Lawvere's functorial semantics, the structuralist pursuit of categoricity turns into an absurdity. For the purpose of this construction is to produce a workable category of models rather than just one model up to isomorphism. In the functorial setting  a theory determines a category, not a structure. This makes the structuralist thinking behind the axiomatic method as expressed by Hilbert in the above quote (Section 3) irrelevant. In the new setting:} 
\begin {quote} The theory appears itself as a generic model \cite{Lawvere:1963-2004}.\end {quote} 
\paragraph{This means that the older structuralist distinction between abstract ``formal" axiomatic theories, on the one hand, and their semantics, on the other hand, doesn't apply; what distinguishes a theory form its (other) models is its generic character rather than formal or abstract character.}  

\paragraph{The setting of Sketch theory is similar to that of Lawvere's Functorial semantics but in the former case generic categories are designed as ``generic shapes" or ``generic figures" rather than axiomatic theories. Unlike the case of Functorial semantics such generic categories are not supposed to have logical properties; in some approaches sketches are not even categories but directed graphs with an additional structure. It seems natural to think of sketches as ``proto-structures" but this is somewhat misleading insofar as the usual notion of structure is concerned. A  sketch doesn't represent a bunch of isomorphic systems but generates non-isomorphic systems (its models). These systems share their generic shape not  in the same sense in which different systems are said to share the same mathematical structure.  In fact they share a shape in a more straightforward sense: a given sketch is a common source of all of its models (i.e. specific functors from this given sketch to a background category). To ``have the same source" is obviously an equivalence relation but this equivalence relation doesn't support anything like the structural abstraction. Unlike a shared structure a shared sketch is  concrete (it is usually even supposed to be finite and easily pictured) while things generated by a sketch can be indeed described as abstract structures in the older sense because they are usually distinguished only up to isomorphism! Thus Sketch theory turns Structuralism upside down and in certain aspects reminds of more traditional ways of doing mathematics. Euclid's geometrical universe is generated by two generic figures, namely, the straight line and the circle, which is tantamount to saying that every geometrical object is constructed  by ruler and compass. The analogy seems to me straightforward.}
\footnote{Does this mean that Ehresmann misconceived of his own invention when he thought of Sketch theory as a general theory of structure? I don't think so. A general theory of structure should not be necessarily a structural theory and should not provide a support for Structuralism as a philosophical view about mathematics.}   

\paragraph{Whether or not the new categorical approach to theory-building  -  differently realized in Functorial semantics, Sketch theory, and the \emph{basic theory} of \cite{Lawvere:1966} - can compete with the standard Hilbert-style structural approach remains an open question. The considered constructions don't allow one to claim that this new approach can work independently: we have seen that Lawvere's  \emph{basic theory} depends on the structural \emph{elementary theory}, Functorial semantics is developed by this author similarly in two steps, and Sketch theory in its existing form uses Set theory and usually doesn't make foundational claims at all. However there is no reason either to claim that the   pre-theoretical notion of collection involved in the standard set-theoretic foundations is indispensible in foundations of mathematics (cf. Mayberry's argument in Section 9 above). It can be replaced by a primitive pre-theoretic notion of category that involves common intuitions about processes (transformations) and their composition. What remains a problem is how to upgrade this pre-theoretic notion of category to a theoretical one without using other means but properly categorical.} 
\paragraph{Which means and constructions may qualify as ``properly categorical'' in a foundational context also remains an open question but I think that the standard machinary of first-order logic used in \cite{Lawvere:1966} and later in \cite{McLarty:1991} for writing down axioms of Category theory after the example of Set theory does not qualify as such. Category theory suggests a change of the traditional conception of logic, which is analogous to the change of the traditional conception of geometry that occurred in the 19th century when people stopped thinking about ``the" geometrical space as a universal container of geometrical objects and learned to think about spaces as objects and about objects as spaces with the notion of \emph{intrinsic geometry} of a given geometrical object. In the first half of the 20th century people learned to think about systems of logic as objects living in larger meta-logical frameworks. Category theory showed how one can think about objects (i.e., appropriate categories) as systems of logic with the notion of internal language of a given category \cite{Lambek&Scott:1986}. This reciprocal move that allows one to avoid the bad infinity of meta-meta.....-logics and meta-meta....-mathematics in foundations of mathematics has immense philosophical importance and I think that it has to be taken into account in categorical foundations.  This is why the presence of a self-standing system of logic representing alleged universal laws of reasoning seems me inappropriate in categorical foundations. A candidate for replacement can be a version of Sketch logic developed in the vein of \cite{Wells&Bagchi:forthcoming} and  \cite{Diskin&Wolter:forthcoming}. I leave this issue for a further study. }

\section{Conclusion: a categorical perspective in and on mathematics}

\paragraph{I hope to have convinced the reader that the project of categorical foundations requires a new philosophical view on mathematics, which the traditional Structuralism cannot possibly provide. Let me now try to summarize this new categorical view by contrasting it with the structuralist view.  What matters in the categorical mathematics is how mathematical objects and constructions transform into each other, not what (if anything) remains invariant under these transformations.  So categorical mathematics is a theory of abstract transformation, not a theory of abstract form. A  theory in categorical mathematics is a generic model (Lawvere) rather than a scheme (Hilbert).} 

\paragraph{The categorical view on mathematics - as distinguished from categorical foundations of mathematics in the sense articulated in the previous Section - suggests a new understanding of the role of history of mathematics in mathematics itself. In the above Introduction I argued that older and recent formulations of the Pythagorean theorem cannot be identified through a merely linguistic translation. Here I strengthen this claim as follows: different versions of this theorem don't share any invariant content; such invariant content is not necessary in order to qualify these formulations as formulations of the same theorem. For the reader's convenience I quote here this example again:}

\begin{quote} \textbf{(1)} In right-angled triangles the square on the side subtending the right angle is equal to the squares on the sides containing the right angle. \end{quote}
\paragraph{(Proposition 1.47 of Euclid's \emph{Elements})} 

\begin{quote} \textbf{(2)}   If two non-zero vectors $x$ and $y$ are orthogonal then$(y - x)^2 = y^2 + x^2$ \end{quote}

\paragraph{(\cite{Doneddu:1965}, slightly modified as explained in Footnote 2):}

\paragraph{Since there is no linguistic translation between \textbf{(1)} and\textbf{(2)}, one may look for a translation of a different sort. Such translations certainly exist but (and this is crucial for my argument) they all lead in one direction, namely from \textbf{(1)} to \textbf{(2)}! There are many ways in which Euclid's geometry can be interpreted in modern terms; \cite{Doneddu:1965} is just one way of doing this among many others. But there is evidently no way to spell modern geometry in Euclid's terms. (As I have explained in the Introduction we must talk here about translation between corresponding theories, not only about translation between separate propositions.) Our history in general and our intellectual history in particular develops from the past to the future, not the other way round; there is no symmetry between the past and the future. Since sound translations between \textbf{(1)} and \textbf{(2)} go in one direction none of them is reversible. According to the argument given in the Section 6, this implies that no such translation allows for the identification of an invariant. Thus the existence of sound translations between theories doesn't imply that these theories share anything like an invariant content.  There is no essence, no conceptual core preserved by the translation of \textbf{(1)} into \textbf{(2)}. But why in this case should we count them as different versions of the same theorem?}

\paragraph{My answer is this: the Pythagorean theorem (as distinguished from its formulations by Euclid and Doneddu) is a particular component of the translation of Euclid's geometry into modern terms proposed by Doneddu, which takes  \textbf{(1)}  into  \textbf{(2)}. I have in mind the categorical notion of the component of a functor explained above in the Section 9.  A similar criterion of identity can be applied, \emph{mutatis mutandis}, to theories and to particular mathematical objects. As far as we are talking about translations like \cite{Doneddu:1965}, which are made with the intent to modernize older mathematical contents, my proposed understanding of the identity of this content amounts, roughly, to the identification of this content with its conceptual history. From this point of view different versions of the Pythagorean theorem produced at different times with different foundations can be described as different temporal stages of the same conceptual entity (i.e. the same collective intellectual process) persisting through time. Noticeably, this persisting entity doesn't reduce to the set of its temporal stages because such a reduction forgets about translations between these stages, which from my categorical viewpoint are the most essential. A mathematical notion is not a set of its temporal stages but a category of translations between these stages. I would like also to stress that in spite of the tendency to consensus alternative foundations of mathematics and alternative theories based on these different foundations can be well developed by contemporaries. So what has been just said about translations between mathematical theories produced at different epochs also applies to alternative contemporary approaches.} 

\paragraph{The reader may wonder if by promoting categorical foundations of mathematics and insisting on the importance of the pedagogical aspect of foundations (in the above Introduction) I am suggesting another radical reform of school mathematics similar to that suggested by enthusiasts  of set-theoretic structural mathematics in the 1960s. This  latter reform, undertaken simultaneously in many countries and known in the US under the name of \emph{New Maths}, was an attempt to conform the elementary mathematical education to the current views on foundations: kids were supposed to learn about sets and elementary structures rather than solve traditional geometrical problems with compasses and ruler and the like \cite{Adler:1972}. Later this reform was commonly recognized as a pedagogical failure \cite{Kline:1973}.}

\paragraph{Although basic elements of Category theory just like elements of Set theory can be learned at early stages of mathematical education (see \cite{Lawvere&Schanuel:1997}) the categorical view on mathematics outlined above suggests a different pedagogical strategy. An average today's elementary mathematical textbook presents a mixture of patterns of mathematical reasoning coming from various historical epochs: it usually contains some elements of Euclid-style synthetic geometry, some Cartesian-style elementary symbolic algebra, some oldish analytic geometry, some elements of Set theory. In my view such a diversity of styles and contents is quite justified because older approaches continue play a role in today's mathematical practice. Pupils definitely need to learn something about older ways of doing mathematics, not just learn older mathematical contents presented in a recent fashionable form. In other words, today's foundations of mathematics must have a historical dimension and include most of the significant older foundations. Recall Lawvere's words quoted in the Introduction concerning the role of foundations in clarifying the ``origins and generals laws of development" of a science.} 

\paragraph{A major problem with such a historical approach in mathematics education is that it appears to be incompatible with the usual notion of axiomatic method. Facing this problem, today's textbooks, written after the failure of the New Maths, often compromise severely against the classical standard of systematicity without suggesting any replacement.}     

\paragraph{My proposed solution of this problem is this: a modern mathematical textbook should provide a few different versions of axiomatic method rather than one. Basic patterns of Euclid-style, Cartesian-style, Hilbert-style and Bourbaki-style theory-building should be definitely included. Radical differences between these approaches should be articulated rather than hidden. The controversial dialectical nature of foundations of mathematics should be stressed rather then kept as a Pythagorean secret. Categorical foundations should be treated not only as a basis of an important part of recent mathematics but also as a means of organization the above diverse contents into a systematic whole. Teaching different ways of doing mathematics the teacher should show how to translate between mathematical contents developed on different foundations without trying to smooth the differences between these foundations.  Given today's overwhelming flow of rapidly updating information such a translation skill is crucially needed in any domain of activity including mathematics. A coherent translation is possible even when no invariant structure is available.} 

\bibliographystyle{plain} 
\bibliography{catstruct} 

\begin{thebibliography}{10}

\bibitem{Adler:1972}
I.~Adler.
\newblock {\em The New Mathematics}.
\newblock John Day and Co (New York), 1972.

\bibitem{Awodey:1996}
S.~Awodey.
\newblock Structure in mathematics and logic: A categorical perspective.
\newblock {\em Philosophia Mathematica}, 3(4):209 -- 237, 1996.

\bibitem{Bourbaki:1939-1983}
N.~Bourbaki.
\newblock {\em Elements de mathematique}.
\newblock Masson (Paris), 1939-1983.

\bibitem{Bourbaki:1950}
N.~Bourbaki.
\newblock The architecture of mathematics.
\newblock {\em The American Mathematical Monthly}, 57(4):221 -- 232, 1950.

\bibitem{Diskin&Wolter:forthcoming}
Z.~Diskin and U.~Wolter.
\newblock A diagrammatic logic for object-oriented visual modeling,
  forthcoming.
\newblock To appear in ENTCS.

\bibitem{Doneddu:1965}
A.~Doneddu.
\newblock {\em Geometrie Euclidienne}.
\newblock Plane (Paris), 1965.

\bibitem{Euclides:1883-1886}
Euclides.
\newblock {\em Opera omnia}.
\newblock Lipsiae, 1883-1886.

\bibitem{Frege:1884}
G.~Frege.
\newblock {\em Die Grundlagen der Arithmetik, eine logisch mathematische
  Untersuchung Ÿber den Begriff der Zahl}.
\newblock Breslau: W. Koebner, 1884.

\bibitem{Frege:1971}
G.~Frege.
\newblock {\em On the Foundations of Geometry and Formal Theories of
  Arithmetic}.
\newblock Yale University Press, 1971.

\bibitem{Hellman:forthcoming}
G.~Hellman.
\newblock Structuralism, mathematical, 1978.
\newblock forthcoming in The Encyclopedia of Philosophy 2d Ed. (MacMillan),
  available at www.tc.umn.edu/~hellm001/.

\bibitem{Hellman:2001}
G.~Hellman.
\newblock Three varieties of mathematical structuralism.
\newblock {\em Philosophia Mathematica}, 2(9):184--211, 2001.

\bibitem{Hilbert:1980}
D.~Hilbert.
\newblock {\em Foundations of Geometry}.
\newblock OpenCourt, 1980 (1899).

\bibitem{Isbell:1967}
J.R. Isbell.
\newblock Review of lawvere:1966.
\newblock {\em Mathematical Reviews}, 34, 1967.

\bibitem{Kline:1973}
M.~Kline.
\newblock {\em Why Johnny CanÕt Add: the Failure of New Maths}.
\newblock St James Press (New York, London), 1973.

\bibitem{Lambek&Scott:1986}
J.~Lambek and Scott~P. J.
\newblock {\em Introduction of Higher-Order Categorical Logic}.
\newblock Cambridge University Press, 1986.

\bibitem{Lawvere:1963-2004}
F.W. Lawvere.
\newblock Functorial semantics of algebraic theories, ph.d. thesis (1963),
  columbia university, 1963-2004.
\newblock republished in Reprints in Theory and Applications of Categories, No.
  5 (2004):1-121, with a new authorÕs Commentary; available at
  www.tac.mta.ca/tac/reprints/articles/5/tr5.pdf.

\bibitem{Lawvere:1964}
F.W. Lawvere.
\newblock Elementary theory of the category of sets.
\newblock {\em Proceedings of the National Academy of Science},
  52(6):1506--1511, 1964.

\bibitem{Lawvere:1966}
F.W. Lawvere.
\newblock The category of categories as a foundation for mathematics.
\newblock In {\em Proceedings of the La Jolla Conference on Categorical
  Algebra}, pages 1--21, 1966.

\bibitem{Lawvere:2003}
F.W. Lawvere.
\newblock Foundations and applications: Axiomatization and education.
\newblock {\em The Bulletin of Symbolic Logic}, 9(2):213--224, 2003.

\bibitem{Lawvere&Schanuel:1997}
F.W. Lawvere and S.H. Schanuel.
\newblock {\em Conceptual Mathematics: a first introduction to categories}.
\newblock Cambridge University Press, 1997.

\bibitem{MacLane:1996}
S.~MacLane.
\newblock Structure in mathematics.
\newblock {\em Philosophia Mathematica}, 4(2):174--183, 1996.

\bibitem{Mayberry:2000}
J.P. Mayberry.
\newblock {\em The foundations of mathematics in the theory of sets}.
\newblock Cambridge University Press, 2000.

\bibitem{McLarty:1991}
C.~McLarty.
\newblock Axiomatizing a category of categories.
\newblock {\em The Journal of Symbolic Logic}, 56(4), 1991.

\bibitem{Popper:1978}
K.~Popper.
\newblock Three worlds, 1978.
\newblock www.tannerlectures.utah.edu/lectures/.

\bibitem{Rodin:2007}
A.~Rodin.
\newblock Identity and categorification.
\newblock {\em Philosophia Scientiae}, 11(2):27--65, 2007.

\bibitem{Rodin:forthcoming}
A.~Rodin.
\newblock How mathematical concepts get their bodies, 2009.
\newblock forthcoming in Topoi.

\bibitem{Wells:1993}
Ch. Wells.
\newblock Sketches: Outline with references, 2009.
\newblock http://www.cwru.edu/artsci/math/wells/pub/pdf/sketch.pdf.

\bibitem{Wells&Bagchi:forthcoming}
Ch. Wells and A.~Bagchi.
\newblock Graph-based logic and sketches, forthcoming.
\newblock a book in progress available at
  http://www.cwru.edu/artsci/math/wells/pub/pdf/gbls.pdf.

\end{thebibliography}
\end{document}